\documentclass[onecolumn]{article}
\usepackage[utf8]{inputenc}
\usepackage[english]{babel}
\usepackage{amsthm}  
\usepackage{amsmath} 
\usepackage{amssymb} 
\usepackage{mathabx} 
\usepackage{graphicx}
\usepackage{fancyhdr}
\usepackage{setspace}
\usepackage{lineno}	 
\usepackage{lastpage}
\usepackage{enumitem}
\usepackage{mathtools}
\usepackage{epsfig} 
\usepackage[outdir=./]{epstopdf} 
\usepackage[margin=1in]{geometry}
\usepackage{caption}  
\usepackage[numbers,comma,sort&compress,super]{natbib}
\usepackage{tocloft}	
\usepackage{hyperref}
\usepackage{cleveref}

\newcommand{\mytitle}{Optimization with Trained Machine Learning Models Embedded}
\newcommand{\myshorttitle}{Optimization with Trained Machine Learning Models Embedded}
\newcommand{\myauthor}{Artur M. Schweidtmann$^{1}$, Dominik Bongartz$^{2,3}$, Alexander Mitsos$^{2,*}$} 
\newcommand{\myauthorshort}{Schweidtmann, Bongartz, Mitsos}
\author{\myauthor}

\allowdisplaybreaks

{
	\fancyhead[L]{\myshorttitle}
	\fancyhead[R]{\small{\textit{\the\day.\the\month.\the\year}}}
	\fancyfoot[L]{\copyright \small{\textit{\myauthorshort}}}
	\fancyfoot[R]{Page \thepage\ of \pageref*{LastPage}}
}

\fancypagestyle{firststyle}
{
	 \fancyhead[L]{\copyright \small{\textit{\myauthorshort}}}
   \fancyhead[R]{Page \thepage\ of \pageref*{LastPage}}
	 \fancyfoot[L]{\footnotesize \rule{0.25\textwidth}{0.4pt} \newline Corresponding author: $^*$A. Mitsos, E-mail: amitsos@alum.mit.edu}
   \fancyfoot[C]{}
	 \fancyfoot[R]{}
}

\newcommand{\E}{\mbox{I\negthinspace E}}

\renewcommand{\vec}[1]{\mathbf{#1}}
\newcommand{\mpmin}[1]{\underset{#1}{\text{min}}~}

\begin{document}

\thispagestyle{firststyle}
\begin{flushleft}\begin{large}\textbf{\mytitle}\end{large} \end{flushleft}
\myauthor

\begin{flushleft}\begin{small}
$^1$ Delft University of Technology, Department of Chemical Engineering, Van der Maasweg 9, Delft 2629 HZ, The Netherlands\\

$^2$ Process Systems Engineering (AVT.SVT), RWTH Aachen University, Forckenbeckstr. 51, 52074 Aachen, Germany. \\

$^3$ Department of Chemical Engineering, KU Leuven, Celestijnenlaan 200F, 3001 Leuven, Belgium \\

\end{small}
\end{flushleft}

\thispagestyle{firststyle}



\section{Introduction}
\label{sec:Introduction}

Machine learning (ML) gives computers the ability to learn from data without being explicitly programmed. 
Cases where the data includes input vectors (in $\mathbb{R}^{n_d}$) and corresponding target vectors (in $\mathbb{R}^{n_t}$) are known as supervised learning.\cite{Bishop.2009} 
Regression problems learn continuous variables as targets whereas classification problems learn discrete categories as targets. 
While linear regression goes back to the least-squares method by Gauss and Legendre in the early 19th century\cite{stigler1981gauss}, recent advances in computing led to the breakthrough performance of ML in various domains even surpassing human performance.\cite{lecun2015deep}

Across many engineering and science domains, ML models are trained on data and are subsequently optimized. 
In chemistry, for instance, ML models are commonly fitted to experimental data. 
This training usually constitutes an optimization problem where model parameters are varied to minimize the model error on the training data. 
Then, a subsequent optimization problem can be solved to identify the experimental conditions (i.e., inputs of the ML model) that maximize the experimental performance (i.e., the output of the ML model). 
We refer to this subsequent optimization problem as \textit{optimization with trained ML models embedded}. 
Notably, the ML model parameters are optimization variables during the training and fixed parameters during the \textit{optimization with trained ML models embedded}. 
Vice versa, the inputs of the ML model are fixed parameters (determined by the training data) during the training and optimization variables during the \textit{optimization with trained ML models embedded}. 
While other chapter in this book deal with the training of ML models (c.f. Section~\ref{sec:see_also}), this chapter focuses on the (deterministic global) \textit{optimization with trained ML models embedded}.  
Also, this chapter covers the optimization of hybrid models\cite{von2014hybrid} where ML models are combined with mechanistic models. 

\section{Formulations}
\label{sec:Formulations}
\begin{figure*}
    \centering
    \includegraphics[width=\textwidth]{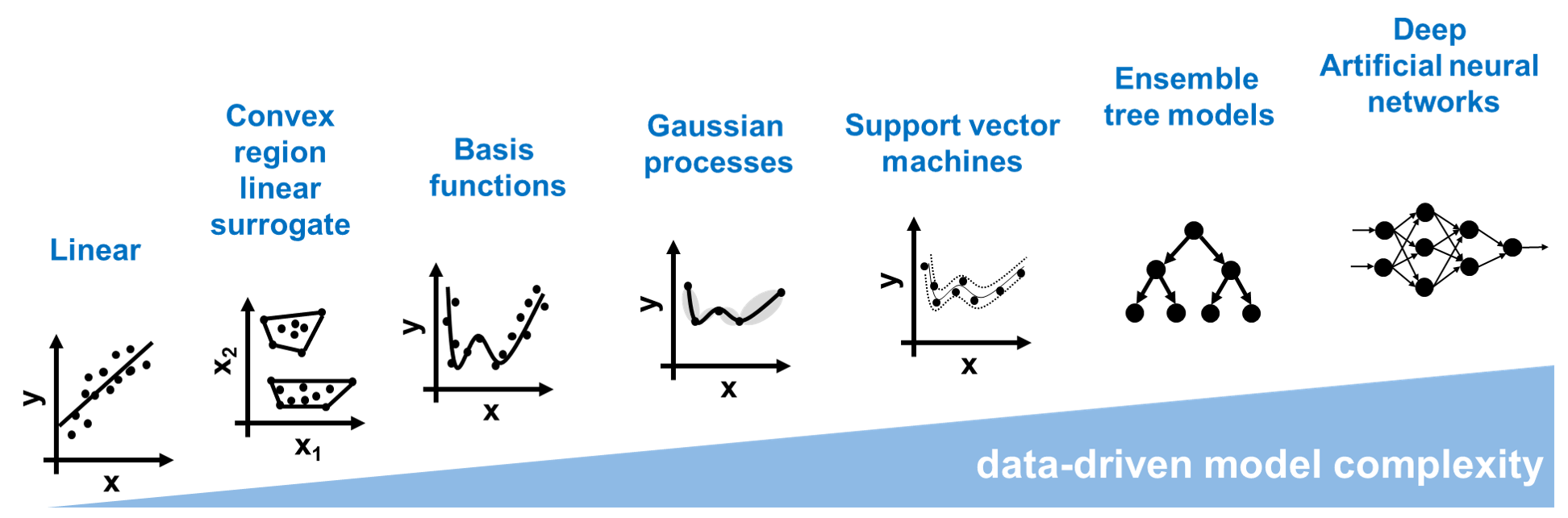}
    \caption{Overview of ML models ordered by increasing complexity. Repreinted from Schweidtmann et al. (2021) \cite{schweidtmann2021machine} (Figure licensed under Creative Commons CC BY).}
    \label{fig:Complexity_of_ML_models}
\end{figure*}

As illustrated in Figure~\ref{fig:Complexity_of_ML_models}, there exists a broad variety of ML models ranging from linear regression models to deep artificial neural networks (ANNs). 

Once a ML regression model has been trained (and validated) on input data from $\mathcal{D}\subset \mathbb{R}^{n_d}$ and target data from $\mathbb{R}^{n_t}$, the trained model usually represents an explicit function $\vec{m}: \mathcal{D}\to\mathbb{R}^{n_t}$. 
Depending on the chosen model, this function may have different properties, e.g., it may be
\begin{itemize}
    \item discontinuous (e.g., for ensemble tree models) or continuous (e.g., ANNs),
    \item nonsmooth (e.g., for ANNs with ReLU activation function) or smooth (e.g., for ANNs with hyperbolic tangent activation function),
    \item linear or convex nonlinear or nonconvex.
\end{itemize}

\subsection{Reduced-Space and Full-Space Model Formulations}
\label{sec:Formulations_RS_FS}
For most ML models, there are multiple approaches of embedding the function $\vec{m}$ in the optimization problem to be solved. 
One way to classify these approaches is to distinguish between \emph{reduced-space formulations} and \emph{full-space formulations}.\cite{bongartz2019deterministic} 

In reduced-space formulations, the ML model is used directly as the function $\vec{m}$.
For example, if the problem to be solved consists in minimizing the output of the ML model over a compact input set $\mathcal{X}\subset\mathcal{D}$, the corresponding reduced-space formulation would read
\begin{equation*}
\mpmin{\vec{x}\in\mathcal{X}} \vec{m}(\vec{x}).
\end{equation*}
Thus, the only variables are the $n_d$ inputs of the model, and there are no constraints beyond those imposed by the set $\mathcal{X}$.
If the ML model is combined with other functions, e.g., in a hybrid model\cite{von2014hybrid} containing both ML and mechanistic parts, the problem could take the form
\begin{align*}
\mpmin{\vec{x}\in\mathcal{X}, \vec{y}\in\mathcal{Y}} &f(\vec{x},\vec{y},\vec{m}(\vec{x}))\\
\text{s.t.} &\vec{h}(\vec{x},\vec{y},\vec{m}(\vec{x}))=\vec{0}\\
&\vec{g}(\vec{x},\vec{y},\vec{m}(\vec{x}))\leq\vec{0},
\end{align*}
where $\vec{y}\in\mathcal{Y}\subset\mathbb{R}^{n_y}$ are additional optimization variables, and $f$, $\vec{h}$, and $\vec{g}$ are some objective function, equality constraints, and inequality constraints, respectively, all of which may depend on the outputs of the ML model as well as the variables $\vec{x}$ and $\vec{y}$. 
Alternatively, the inputs of the ML model may not be optimization variables themselves, but rather outputs of some other function $\vec{i}: \mathcal{Y}\to\mathcal{X}$:
\begin{align*}
\mpmin{\vec{y}\in\mathcal{Y}} &f(\vec{y},\vec{m}(\vec{i}(\vec{y})))\\
\text{s.t.} &\vec{h}(\vec{y},\vec{m}(\vec{i}(\vec{y})))=\vec{0}\\
&\vec{g}(\vec{y},\vec{m}(\vec{i}(\vec{y})))\leq\vec{0}.
\end{align*}

In a full-space formulation, internal variables of the ML models become part of the optimization problem as optimization variables, and equations defining the ML model become constraints. In this case, the function $\vec{m}$ does not occur in the problem as a single function that can be evaluated to obtain outputs for given inputs, but instead, the evaluation becomes part of solving the optimization problem.
Note that typically  a range of full-space formulations exists that differ in how many and which internal variables and equations become part of the optimization problem. 
These can also be considered as intermediate versions between full- and reduced-space formulations. 
Similar to the reduced-space formulation, more complex problems can of course be constructed where the ANN is combined with other functions, e.g., in a hybrid mechanistic/data-driven model.\cite{von2014hybrid}

From a modeling perspective, reduced-space formulations are very intuitive and easy to implement in procedural environments such as e.g., when modeling in Python. In these environments, the ML models are naturally implemented as functions that can just be used in the optimization problem. In contrast, when modeling in an equation-oriented environment such as, e.g., GAMS or AMPL, the reduced-space formulation is often much more challenging to implement. For example, to write an ANN in a reduced-space formulation, one would have to symbolically substitute the outputs of each layer into the activation function of the next layer in order to ultimately obtain an expression for the outputs of the ANN (i.e., those of the last layer) as a function of its inputs (i.e., the inputs of the first layer). In such an environment, it is more natural to write a full-space formulation that merely requires defining the outputs of all layers as variables and writing the corresponding equations as constraints.

From a computational perspective, the chosen formulation can have a considerable impact on the effort for solving the problem, because the problem looks very different to the optimizer: The reduced-space formulations have much fewer variables and equality constraints than the full-space formulations but contain potentially very complicated functions.
In particular, the number of variables and constraints in the optimization problem is independent of the type of ML model used or the amount of data used for training it. In the full-space formulation, the number of optimization variables and constraints may increase with the complexity of the ML model (e.g., for an ANN, it may increase with the number of hidden layers, or with the number of neurons per layer\cite{Schweidtmann2019detGlobalANN}), or even with the amount of data that was used for training the model (e.g., for Gaussian processes\cite{schweidtmann2021deterministic}).
As discussed in the following, it may depend on the optimization algorithm how efficiently either formulation can be solved.

An alternative classification of ways of incorporating ML models in optimization problems is given by Lombardi \& Milano\cite{lombardi2018boosting} who distinguish \emph{embedding an ML model by expanding the language} as opposed to \emph{encoding the ML model using the native language}. 
This first case entails defining a new function representing the ML model that can then be used when formulating the optimization problem. If this is done for an entire ML model, this would correspond to a reduced-space formulation as defined above. 
This approach may, however, only introduce certain components of an ML model as new functions, e.g., activation functions of ANNs. 
The second case entails representing the ML model via functions that are already available in the chosen optimization environment. This could in fact be done for both full space and reduced space formulations. 
In full space formulations, the ML model consists of a number of variables and constraints in the optimization problems, while in RS formulations, the functions encoding the ML model are a composition of existing functions. 

\subsection{Modeling Nonsmooth Functions}
\label{sec:Formulations_Nonsmooth}
For some ML models such as ANNs with ReLU activation function, there is an additional distinction regarding the way the nonsmoothness is handled. 
On the one hand, the ReLU activation function can be formulated with the $\max$ function. 
This enables a reduced-space formulation of the ANN, meaning that the network can be expressed as an explicit function mapping the inputs to the outputs of the network. 
However, this function is usually nonconvex and nonsmooth and thus challenging to solve to global optimality. 
On the other hand, the $\max$ function can be reformulated as a piecewise linear function modeled with a binary variable that determines which of the linear pieces is active (c.f. Section~\ref{sec:neural_network}).
This avoids the use of nonsmooth functions in the model at the expense of adding one binary variable per neuron with ReLU activation function. 
Furthermore, this makes a reduced-space formulation as described above impossible, because the optimizer now needs to control the binary variables, which are part of the internals of the ANN. 
Typically, this also entails the use of additional constraints, e.g., in the form of a Big-M formulation.
Similarly, the nonsmoothness of tree models can be reformulated using binary variables (c.f. Section~\ref{sec:decision_tree}). 

\subsection{Validity Domain Modeling}
\label{sec:Formulations_Validity}
An additional aspect is related to the domain of the function $\vec{m}$. 
Since the ML model has been trained on a limited number of data points $(\vec{x}_i,\vec{t}_i) \in \mathbb{R}^{n_d}\times\mathbb{R}^{n_t}$, it can not be expected to give a good description of the underlying, ``true'' input-output relationship for arbitrary inputs $\vec{x}\in\mathbb{R}^{n_d}$. 
If $\vec{x}$ is far from the inputs of the training data, the model needs to extrapolate, which can give unreliable results. 
In this sense, it may be prudent to add constraints to the optimization problem that restrict the inputs of the ML model to values that are not too far from its training data. 
Such constraints could take the form of the convex hull of the training data points as in the convex region surrogate models\cite{Zhang.2016b}, or an approximation of the possibly nonconvex region covered by the training data points, e.g., obtained by one-class classifiers.\cite{schweidtmann2021obey} 
Also, the distance to training data points can be included as a penalty to the objective function.\cite{mistry2020mixed} 

\section{Machine Learning Models}
In this section, we discuss a few common ML models and possible optimization formulations with these ML models embedded. 

\subsection{Convex region surrogate models}
Convex region surrogate models were proposed by Zhang et al.\cite{Zhang.2016b} 
They consist of affine functions defined over a number of convex polytopes that are constructed to contain the training data. They can approximate nonlinear (potentially discontinuous) functions over nonconvex regions. 
Given the regions into which the domain is divided according to the model, evaluating the model at a point amounts to identifying which region this point lies in and evaluating the affine function defined on this region.

The models are originally devised in a Generalized Disjunctive Programming~(GDP) setting, but they can be reformulated as MILPs, e.g., via the hull reformulation. 
This representation is particularly suitable if the optimization problem in which the model is to be embedded is already a MILP.
According to the above classification, this is a full-space formulation of the model, because the optimization solver controls additional variables and constraints.

To our knowledge, a reduced-space formulation of this model has not been attempted yet. In its purest form, this would entail expressing the model as a piecewise-defined function of $\vec{x}$. 
In this case, the challenge lies in the fact that the function may be discontinuous, and that the constraints on the input values need to be accounted for explicitly.

\subsection{Gaussian processes}
\label{sec:gaussian_process}
Gaussian processes, also known are Kriging, are infinite-dimensional generalizations of multivariate Gaussian distributions.\cite{rasmussen2004gaussian}
Their predictions follow Gaussian distributions and thus include an estimate and a variance. 
Formally, Gaussian processes are defined by its mean function ($m(\boldsymbol{x})\coloneqq \E\big[\tilde{f}(\boldsymbol{x})\big]$) and covariance function ($k(\boldsymbol{x},\boldsymbol{x}')\coloneqq \E \big[~( \medspace y(\boldsymbol{x}) - m(\boldsymbol{x}) \medspace )~( \medspace y(\boldsymbol{x}') - m(\boldsymbol{x}) \medspace )^{\mathrm{T}} \big]$), where $y$ is an observation from $\tilde{f}(\boldsymbol{x})$. 
The most common covariance function is the squared exponential covariance function, $k_{SE}'(r) \coloneqq \exp \left( -\frac{1}{2}~r^2 \right)$, where $r \coloneqq \sqrt{(\boldsymbol{x} - \boldsymbol{x}')^{\mathrm{T}} \boldsymbol{\Lambda}~(\boldsymbol{x} - \boldsymbol{x}')}$ is a Euclidean distance with weighting factors $\boldsymbol{\Lambda} \coloneqq \mathrm{diag}(\lambda_1^2, \cdots , \lambda_i^2, \cdots \lambda_{n_x}^2)$. 
The hyperparameters $\boldsymbol{\theta} = [\lambda_1,..,\lambda_d]$ are adjusted during training to maximize the probability that the Gaussian process fits the training data. 
The Gaussian process posterior distribution is obtained by conditioning the prior distributions on training data ($\boldsymbol{\mathcal{X}},\boldsymbol{\mathcal{Y}})$, $\tilde{f}(\boldsymbol{x}) \sim \mathcal{GP}(m(\boldsymbol{x}), k(\boldsymbol{x},\boldsymbol{x}') | \boldsymbol{\mathcal{X}},\boldsymbol{\mathcal{Y}})$. 
The estimate and variance of the posterior can be calculated by
\begin{align}
	m_{\boldsymbol{\mathcal{D}}}( \boldsymbol{x}) &= \boldsymbol{K}_{\boldsymbol{x},\mathcal{X}} \left( \boldsymbol{K}_{\mathcal{X},\mathcal{X}} \right)^{-1}  \boldsymbol{y},  \label{eq:global_gp_optimization_GP_prediction}\\
	k_{\boldsymbol{\mathcal{D}}}(\boldsymbol{x}) &= {K}_{\boldsymbol{x},\boldsymbol{x}} - \boldsymbol{K}_{\boldsymbol{x},\mathcal{X}} \left( \boldsymbol{K}_{\mathcal{X},\mathcal{X}} \right)^{-1} \boldsymbol{K}_{ \mathcal{X},\boldsymbol{x} }, \label{eq:global_gp_optimization_GP_variance}
\end{align}
where $\boldsymbol{K}_{\mathcal{X},\mathcal{X}} \coloneqq \left[ k(\boldsymbol{x}_i,\boldsymbol{x}_j) \right] \in \mathbb{R}^{N \times N}$ is the covariance matrix of the data, 
$\boldsymbol{K}_{\boldsymbol{x},\mathcal{X}} \coloneqq \left[ k(\boldsymbol{x},\boldsymbol{x}_1^{(\boldsymbol{\mathcal{D}})}), ..,k(\boldsymbol{x},\boldsymbol{x}_N^{(\boldsymbol{\mathcal{D}})}) \right] \in \mathbb{R}^{1 \times N}$ is the covariance vector between the candidate point $\boldsymbol{x}$ and the training data, 
and ${K}_{\boldsymbol{x},\boldsymbol{x}} \coloneqq k(\boldsymbol{x},\boldsymbol{x})$.

Trained Gaussian processes are commonly embedded in optimization problems for iterative planning of (computational) experiments through Bayesian optimization.
Therein, an acquisition function is maximized to identify an optimal tradeoff between exploration and exploitation. 
One common example is the maximization of the expected improvement acquisition function.\cite{jones1998efficient} 
Furthermore, Gaussian processes are commonly used as surrogate models when only limited data is available.\cite{caballero2008rigorous,boukouvala2017argonaut}
Optimization problems with trained Gaussian processes embedded are NLPs as the covariance function is nonlinear. 
Moreover, the optimization of acquisition functions are known to exhibit a large number of (suboptimal) local optima.\cite{kim2019local}
In the previous literature, optimization problems with trained Gaussian processes embedded have mostly been solved using full space formulations (e.g.\cite{caballero2008rigorous,boukouvala2017argonaut}) while recent work indicated that the reduced space formulation can be advantageous.\cite{schweidtmann2021deterministic}
Notably, tailored convex and concave relaxations of covariance and acquisition functions can improve the solution times for these problems significantly.\cite{schweidtmann2021deterministic}

\subsection{Decision tree ensembles}
\label{sec:decision_tree}
Decision tree ensembles are popular models for representing nonlinear relationships in large-scale, noisy data sets including both continuous and discrete inputs. 
They consist of multiple decision trees, where a decision tree predicts the output based on a sequence of queries involving the inputs (e.g., if $x>c$, where $c$ is a constant parameter, follow branch 1 of the tree, if $x\le c$ follow branch 2 of the tree). 
For given values of the inputs, the output is computed by following the tree based on the answers to those queries from the root to a certain leaf, which then contains the output value. 
In a tree ensemble, the output of multiple trees is averaged (in the case of regression) to yield the output of the overall model to improve the prediction accuracy. 
The most popular models of this type are gradient boosted trees, where the trees are built sequentially by learning the residual between the data and the previous tree, and random forests, which consist of trees trained on random subsets of the data as well as the inputs used in the queries.

By construction, tree ensemble models predict the output as a piecewise constant function of the inputs. 
As such, they are inherently discontinuous, which makes the optimization of problems with trained ensemble models embedded difficult to optimize. 
Recently, Mi\v{s}i\'{c}\cite{mivsic2020optimization} and Biggs et al.\cite{biggs2017optimizing} presented mixed-integer formulations for representing ensemble tree models in optimization problems. 
They encode the results of the queries as binary variables and introduce constraints that enforce the selection of the correct leaf node based on these variables. 
The approach of Biggs et al.\cite{biggs2017optimizing} has the advantage that additional (linear) constraints are easier to incorporate.

Misener and co-workers recently extended these approaches by incorporating an uncertainty measure as a penalty term in the objective. 
This term penalizes large deviation of inputs from the training data, either considering the distance to the nearest (clusters of) data points\cite{thebelt2021entmoot}, or the distance to a subspace of the input space containing the data.\cite{mistry2020mixed} 
The former approach was also used in the opposite sense for application in Bayesian optimization, i.e., to enhance exploration rather than restrict the search to the vicinity of the training data.\cite{thebelt2021entmoot} 

Optimization with tree ensembles models embedded becomes particularly challenging when the ensembles get large (e.g., several thousand trees) or when the trees contain many queries. 
Therefore, dedicated approaches for large models have been proposed. These are based on Benders decomposition \cite{biggs2017optimizing,mivsic2020optimization}, truncation of the tree depth \cite{mivsic2020optimization}, considering subsets of the trees \cite{biggs2017optimizing, mistry2020mixed}, or dedicated B\&B procedures.\cite{mistry2020mixed,thebelt2021entmoot} 

\subsection{Artificial neural networks}
\label{sec:neural_network}
ANNs are one of the most common ML models and they are an essential block of deep learning architectures~(see also, e.g., \cite{Bishop.2009,lecun2015deep}). 
Feed-forward ANNs, also known as multilayer perceptrons, consist of $N$ (hidden) layers including $D^{(k)}$ neurons. 
Generally, ANNs can be viewed as directed acyclic graphs connecting neurons that are organized in layers. 
The first layer corresponds to the inputs of the ANN and the last layer corresponds to its outputs. 
For each layer $k \ge 2$, the output vector $\textbf{z}^{(k)}$ of the layer $k$ is given by
$\textbf{z}^{(k)} = h_{k} \left( \textbf{W}^{(k-1)} \textbf{z}^{(k-1)} + \textbf{b}^{(k-1)} \right)$ 
where $h_{k}(\cdot)$ is an activation function, $\textbf{W}^{(k-1)}$ is a weight matrix and $\textbf{b}^{(k-1)} \in \mathbb{R}^{D^{(k)}}$ a bias vector. Notably, the weights $(\textbf{W}^{(k)})$ and bias vectors $(\textbf{b}^{(k)})$ are known constant parameters in optimization problems with trained ANNs embedded. 
For regression problem, the hyperbolic tangent function, $h_k(x)=\tanh(x)$, has been the most common activation function in ANNs. 
Also, the rectified linear unit (ReLU) activation function, $h_k(x)=\max(x, 0)$, is increasingly common for deep ANNs as it is more effectively trainable.\cite{arora2016understanding} 

Given the nonlinear activation functions, optimization problems with ANNs embedded are generally nonlinear optimization problems (NLPs). 
As discussed in Section~\ref{sec:Formulations_RS_FS}, ANNs can be modeled using full space or reduced space formulations. 
Previous studies indicate that the reduced space formulation has a favorable performance compared to the full space formulation for deterministic global optimization with trained ANNs embedded.\cite{Schweidtmann2019detGlobalANN} 
This is mostly due to the increased problem size of full space formulations and the difficulty to solve large-scale NLPs to global optimality. 

As the ReLU function is a piecewise linear function, it can also be reformulated through binary variables and linear constraints. 
This allows formulating optimization problems with trained ReLU ANNs as mixed-integer linear problems (MILPs).\cite{fischetti2018deep,cheng2017maximum,tjeng2017evaluating} 
As intermediate variables are introduced to the optimization problem, these MILP formulations are essentially full space formulations. 
While most previous works use big-M formulations for ReLU ANNs, recent works also proposed sharp or ideal formulations\cite{anderson2020strong} as well as hybrids\cite{tsay2021partition}. 
Moreover, most previous works employ specialized bound tightening techniques. 
Recently, open-source modeling tools for ReLU ANNs as MILP have been proposed bridging established ML and optimization tools.\cite{reluMIP,bergman2021janos} 

To the best of our knowledge, a comparison between the nonlinear reduced space formulation and the MILP full space formulation has not yet been published in the literature. 
Thus, it still remains unclear what the most efficient problem formation is.  
While the reduced space formation has significantly smaller problem sizes, the MILP formulation is linear. 
Irrespective of the formulation, it appears that the problem size of previous studies has been limited by high computational times. 

\section{Conclusions}
Trained ML models are commonly embedded in optimization problems. 
In many cases, this leads to large-scale NLPs that are difficult to solve to global optimality. 
While ML models frequently lead to large problems, they also exhibit homogeneous structures and repeating patterns (e.g., layers in ANNs). 
Thus, specialized solution strategies can be used for large problem classes. 
Recently, there have been some promising works proposing specialized reformulations using mixed-integer programming or reduced space formulations. 
However, further work is needed to develop more efficient solution approaches and keep up with the rapid development of new ML model architectures.

\bibliographystyle{abbrv}  
\bibliography{Bibs}

\end{document}